# Comparison of Trefftz-Based PINNs and Standard PINNs Focusing on Structure Preservation

Koji KOYAMADA[1]*, Hiroaki OHTANI[2]

[1]Department of Data Science,Osaka-seikeiUniversity
[2]NationalInstitute for Fusion Sciency
*koyamada@g.osaka-seikei.ac.jp

**Abstract.** In this study, we investigate the capability of physics-informed neural networks (PINNs) to preserve global physical structures by comparing standard PINNs with a Trefftz-based PINN (Trefftz-PINN). The target problem is the reproduction of magnetic field-line structures in a helical fusion reactor configuration. Using identical training data sampled from exact solutions, we perform comparisons under matched mean squared error (MSE) levels. Visualization of magnetic field lines reveals that standard PINNs may exhibit structural collapse across magnetic surfaces even when the MSE is sufficiently small, whereas Trefftz-PINNs successfully preserve the global topology of magnetic field lines. Furthermore, the proposed framework is extended to computational fluid dynamics (CFD) problems, where streamline structures of velocity fields are analyzed. Similar tendencies are observed, demonstrating that Trefftz-PINNs provide superior structure preservation compared to standard PINNs. These results indicate that minimizing numerical error alone does not guarantee physical consistency, and that constraining the solution space prior to learning is an effective strategy for physics-consistent surrogate modeling.

## 1. Introduction

### 1.1. Background

Physics-Informed Neural Networks (PINNs) have attracted significant attention as a framework for solving partial differential equations (PDEs) by embedding physical laws into neural network training through PDE residual minimization. PINNs have been successfully applied to a wide range of scientific and engineering problems, in-cluding fluid dynamics, electromagnetics, and heat transfer, owing to their flexibility and mesh-free nature.

However, recent studies and practical experiences have revealed that minimizing PDE residuals does not necessarily guarantee physically correct solutions. In partic-ular, PINNs may converge to solutions that exhibit small residual values while sig-nificantly deviating from the true solution, especially in problems involving complex solution structures or ill-conditioned training dynamics.

### 1.2. Motivation: Residual Hallucination in PINNs

To clarify this issue, we consider a simple steady-state heat conduction problem governed by the Laplace equation, for which an exact analytical solution is available.This problem serves as a minimal yet illustrative example to demonstrate a fundamental limitation of residual-based learning.

Figure 1 compares the exact solution, a purely data-driven neural network trained with supervised data, and a PINN trained primarily by minimizing the PDE residual.Although the PINN successfully reduces the PDE residual to a nearly zero level, the obtained temperature distribution significantly deviates from the exact solution and the data-driven surrogate.

To further quantify the nature of residual hallucination, we conducted a systematic numerical study by varying activation functions, network depth, and the number of neurons per layer. The target problem is the three-dimensional advection–diffusion equation with a known exact solution. Training data were generated by randomly sampling space–time points and using the exact concentration values as supervision for fully connected deep neural networks.

After training, we evaluated not only the mean squared error (MSE) of the concentration field but also the MSE of the second-order spatial derivatives, i.e., the Laplacian of the concentration. Fig.1 shows the error distributions for representative activation functions, including tanh, sin, and swish, across different network architectures.

While the concentration error exhibits clear dependence on activation functions and network configurations, the error of the second-order derivatives remains almost invariant with respect to both network depth and width, regardless of the activation function employed. This observation indicates that higher-order derivative information is not implicitly learned when the network is trained solely on function value data.

Importantly, despite the accurate reconstruction of the concentration field, the failure to reproduce its second-order derivatives provides direct numerical evidence of residual hallucination in PINNs. This result highlights that apparently small residuals or good agreement in function values do not necessarily guarantee physically consistent derivative structures, underscoring a fundamental limitation of residual-based PINNs.

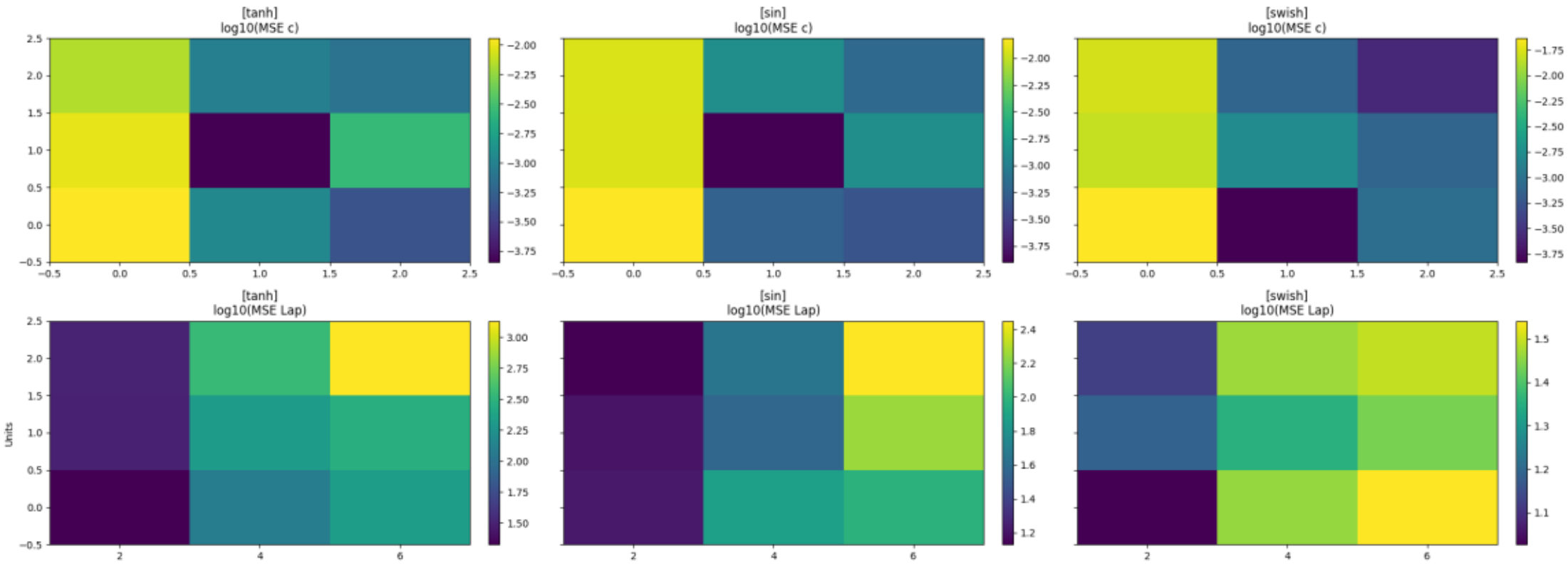

**Figure 1:** Activation Function and Network Architecture Dependence of Second-Order Derivative Errors

This discrepancy is rooted in the theoretical limitations of the Universal Approximation Theorem. While the theorem guarantees that a neural network can converge to a target function (0th-order convergence), it does not guarantee the convergence of its second-order derivatives. Consequently, the PDE residual may diminish prematurely without the network truly capturing the underlying physical curvature, leading to a small residual that fails to enforce consistency with the exact solution.

## 1.3. Implications for Structure-Preserving Problems

While the failure shown in Fig. 2 occurs in a scalar field problem, its implications are more severe in systems characterized by geometric or topological structures.In such problems, including magnetic field line configurations in fusion devices and streamline structures in fluid flows, even small local errors can lead to global structural breakdowns, such as crossing of magnetic surfaces or distortion of flow topology.

Therefore, ensuring physical structure preservation requires not only residual minimization but also appropriate constraints on the admissible solution space.This observation motivates the present study, which investigates the use of **Trefftz bases** to restrict the solution space prior to learning and to suppress residual hallucination.

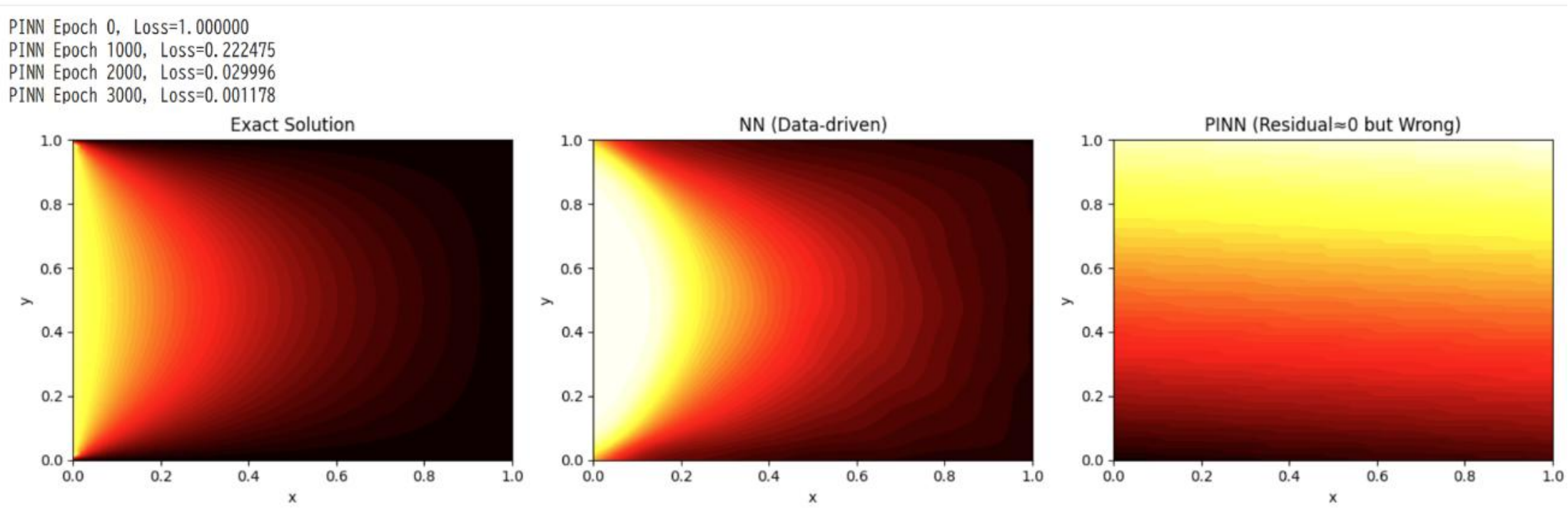

**Figure 2:** Failure case of a physics-informed neural network in a steady-state heat conduction problem governed by the Laplace equation

Although the PDE residual is minimized to a near-zero value, the PINN solution deviates significantly from the exact solution due to gradient vanishing. A purely data-driven neural network trained with supervised data successfully reproduces the correct temperature distribution. This example illustrates the phenomenon of *residual hallucination*.

## 2. Related works

Physics-Informed Neural Networks (PINNs) were originally introduced as a unified framework for incorporating physical laws, expressed as partial differential equa-tions (PDEs), into neural network training through residual minimization [1,2]. Since their introduction, PINNs have been applied to a wide range of forward and inverse problems in computational physics, including fluid dynamics, heat transfer, and electromagnetics.

However, subsequent studies have revealed that PINNs often suffer from severe optimization difficulties. In particular, loss imbalance between PDE residuals and boundary or initial condition terms has been identified as a major obstacle, leading to slow convergence or physically inconsistent solutions [3–5]. Wang et al. analyzed PINN training dynamics from the perspective of neural tangent kernels and showed that minimizing the PDE residual does not necessarily guarantee convergence to the correct solution manifold [3]. These findings indicate that low residual values alone are insufficient to ensure physical correctness.

Another line of research has focused on the spectral bias of neural networks. It has been shown that standard neural networks tend to learn low-frequency components of a solution faster than high-frequency ones, which can result in the suppression of essential spatial structures even when residuals are small [6,7]. This phenomenon is particularly problematic for PDEs whose solutions contain sharp gradients or global geometric structures.

Recent theoretical analyses have shed light on why standard PINNs often fail to converge to physically consistent solutions. Wang et al. analyzed the training dynamics of PINNs from the perspective of the **Neural Tangent Kernel (NTK)**. Their work mathematically demonstrated that a severe discrepancy often exists between the convergence rates of the different loss components (e.g., PDE residuals vs. boundary conditions). This gradient pathology causes the network to prioritize certain loss terms while neglecting others, leading to a state where the PDE residual appears minimized while the global physical structure is entirely misrepresented—a phenomenon we categorize as **"Residual Hallucination."**

While several strategies have been proposed to mitigate these issues, such as self-adaptive weighting of loss functions or respecting causality during training, these methods fundamentally rely on adjusting the optimization process of an unconstrained neural network. In contrast, our proposed **Trefftz-PINN** framework addresses this problem at the architectural level. By restricting the solution space to a physically admissible manifold using Trefftz basis functions a priori, we eliminate the reliance on delicate loss-balancing heuristics. This structural constraint ensures that the model remains robust against the gradient imbalances and spectral biases that typically plague standard residual-based PINNs.

More recently, explicit failure modes of PINNs have been systematically investi-gated. Krishnapriyan et al. demonstrated that PINNs may converge to solutions with negligible PDE residuals while still deviating significantly from the exact solution, especially in diffusion-dominated or stiff problems [8]. Mishra and Molinaro further provided theoretical bounds showing that residual minimization does not directly control the solution error, highlighting a fundamental gap between residual accuracy and solution fidelity [9]. Comprehensive surveys have since summarized these limi-tations and emphasized the need for structural constraints beyond residual-based training [10].

Related to these observations, several studies have proposed solver-in-the-loop or correction-based approaches, in which neural network predictions are explicitly corrected by numerical solvers or projection steps [11,12]. These methods implicitly acknowledge that residual minimization alone may generate unphysical solutions, although they do not explicitly frame this issue as a structural hallucination.

In parallel, reduced-order modeling techniques such as Proper Orthogonal De-composition (POD) have long been used in CFD to restrict the solution space to physically meaningful subspaces derived from snapshots [13,14]. While effective, POD-based models depend heavily on the availability and representativeness of training data and lack explicit enforcement of governing equations.

From a numerical analysis perspective, Trefftz methods offer a fundamentally different approach by constructing solution spaces composed of functions that exactly satisfy the governing PDEs a priori [15–17]. By restricting approximation spaces before optimization, Trefftz methods inherently avoid solutions that violate the un-derlying physics. Recently, hybrid approaches combining Trefftz bases with neural networks have been explored to improve accuracy and stability [18,19].

In contrast to existing PINN failure analyses, the present study introduces the concept of Residual Hallucination, defined as a phenomenon in which neural net-works trained to minimize PDE residuals converge to solutions that are mathemat-ically consistent in a weak residual sense but physically or structurally incorrect. While related to previously reported optimization pathologies [3,8], residual hallu-cination emphasizes the decoupling between residual minimization and global structure preservation, which has not been explicitly highlighted in prior work.

By comparing standard PINNs and Trefftz-PINNs under matched MSE conditions for magnetic field lines and CFD streamlines, this study provides concrete evidence that restricting the solution space before learning is crucial for preventing residual hallucination and preserving physically meaningful structures.

Furthermore, a fundamental limitation arises from the **Universal Approximation Theorem (UAT)**. While UAT guarantees that a neural network can approximate any continuous function to arbitrary precision, it does not inherently ensure the accurate approximation of its **partial derivatives**. Since the loss function of a standard PINN is primarily composed of PDE residuals—which are defined by these derivatives—the lack of theoretical bounds on derivative approximation often leads to "residual hallucination" and non-physical oscillations. This discrepancy between function-value approximation and derivative-consistency is a major factor in the failure of standard PINNs to preserve global structures. To address this, we employ the **Trefftz method** to constrain the solution space. By using basis functions that satisfy the governing equations *a priori*, our approach ensures derivative-consistency by construction, bypassing the inherent limitations of derivative approximation in standard UAT-based architectures.

# 3. Problem Formulation and Methods

## 3.1. Governing Equations

### 3.1.1. Magnetic Field Problem

T We consider a vacuum magnetic field represented by a scalar magnetic potential $\Phi$.The magnetic field $\mathbf{B}$ is given by

$$\mathbf{B} = -\nabla\Phi \tag{1}$$

In current-free regions, the magnetic potential satisfies the Laplace equation

$$\nabla^2\Phi = 0 \tag{2}$$

This formulation is commonly employed in magnetic field analysis of fusion devices, where the plasma region is approximated as a vacuum for field-line analysis. The influence of external coils is incorporated through boundary conditions imposed on $\Phi$, while the coils themselves are assumed to be located outside the computational domain. This assumption allows us to focus on the internal magnetic field structure without explicitly modeling coil currents. To enable quantitative comparison, analytical solutions with helical symmetry that exactly satisfy the Laplace equation are employed as reference solutions.

#### 3.1.2. CFD Problem

For computational fluid dynamics (CFD) analysis, we consider a decaying Taylor–Green vortex as a representative incompressible flow with a known analytical structure. The velocity field $\mathbf{u} = (u, v, w)$is governed by the incompressible Navier–Stokes equations under periodic boundary conditions.

In this study, a single time-step evolution is examined, where the analytical decay rate of the vortex amplitude is known. This setup provides a controlled environment to evaluate whether surrogate models can preserve coherent flow structures, such as vortices and streamlines, beyond pointwise accuracy.

### 3.2. Compared Models

Three types of models are compared in both magnetic field and CFD problems:

1. **Exact solution (True)**: Analytical solutions or reference numerical solutions used as ground truth for validation.
2. **Standard PINNs**: Neural networks trained by minimizing a composite loss function consisting of PDE residuals and data fitting errors. The solution is represented entirely by a neural network without explicit restriction on the admissible solution space.
3. **Trefftz-PINNs**: Hybrid models that combine Trefftz basis functions, which exactly satisfy the governing PDE, with a residual neural networe.

### 3.3. Trefftz-Based Model Construction

In the Trefftz-PINN formulation, the solution $u(\mathbf{x})$is represented as

$$u(\mathbf{x}) = \sum_{i=1}^{N_b} c_i \ \phi_i(\mathbf{x}) + u_{\mathrm{NN}}(\mathbf{x}) \tag{3}$$

where $\phi_i(\mathbf{x})$are Trefftz basis functions that exactly satisfy the governing partial differential equation,$c_i$are scalar coefficients to be optimized, and$u_{\mathrm{NN}}(\mathbf{x})$is a residual neural network responsible for local corrections.

By construction, the dominant component of the solution resides in a physically admissible solution space. The role of the neural network is therefore limited to compensating for boundary effects, discretization artifacts, or modeling errors, rather than representing the global solution structure. This design contrasts with standard PINNs, in which the entire solution is encoded in a neural network and physical constraints are enforced only through loss terms.

### 3.4. Evaluation Strategy

To ensure a fair and meaningful comparison between models, the following conditions are enforced:

- Identical training samples and visualization seeds across all models
- Matched mean squared error (MSE) levels between PINNs and Trefftz-PINNs
- Visualization of physical structures, including magnetic field lines and fluid streamlines

By matching MSE levels, differences in solution quality can be attributed to model design rather than numerical accuracy alone. This evaluation strategy enables us to assess whether constraining the solution space using Trefftz bases improves the preservation of physically meaningful structures in both electromagnetic and fluid flow problems.

## 4. Results

### 4.1. Magnetic Field-Line Structures

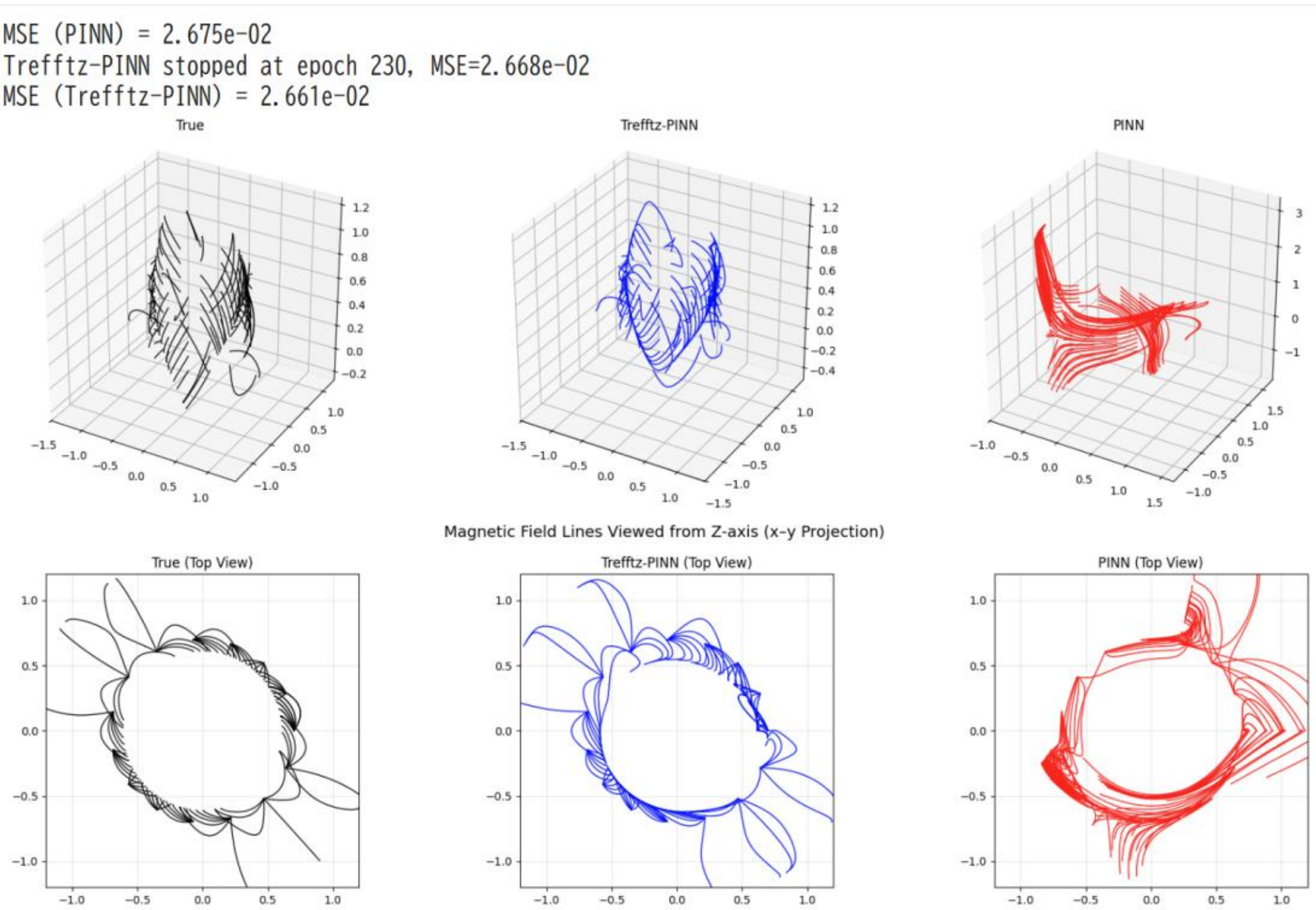

**Figure 3** The magnetic field-line structures reconstructed by the three models: the exact solution, standard PINNs, and Trefftz-PINNs, under conditions where the mean squared error (MSE) is carefully matched between the learning-based models.

**Despite achieving a small MSE comparable to that of Trefftz-PINNs, the standard**

**PINNs exhibit noticeable degradation of magnetic field-line topology. In particular, magnetic field lines generated by PINNs are observed to cross magnetic surfaces, violating a fundamental physical property of vacuum magnetic fields. Such topological inconsistencies occur even though the pointwise numerical error remains small, indicating that MSE alone is insufficient to assess the physical validity of the solution.**

In contrast, Trefftz-PINNs successfully preserve coherent magnetic surfaces that are consistent with the exact solution. Field lines remain confined to their respective magnetic surfaces and exhibit smooth helical structures without spurious crossings. This result demonstrates that constraining the solution space using Trefftz bases effectively suppresses unphysical modes that can arise in standard PINNs.

These observations highlight a critical limitation of residual-based learning: minimizing PDE residuals does not necessarily guarantee the preservation of global field-line topology. Fig. 4 highlights that residual minimization alone does not prevent topological failure, whereas solution-space restriction via Trefftz bases ensures structural fidelity.

## 4.2. MSE versus Number of Trefftz Bases

To investigate the influence of Trefftz basis selection, we examine the relationship between the MSE and the number of Trefftz bases $N_b$. shows that this relationship is non-monotonic (Fig.4).

As $N_b$ increases from a small value, the MSE initially decreases, reflecting the enhanced expressiveness of the Trefftz representation. However, beyond an optimal range, further increasing the number of bases leads to stagnation or even degradation of accuracy. This behavior is attributed to over-parameterization and adverse interactions between the Trefftz expansion and the residual neural network.

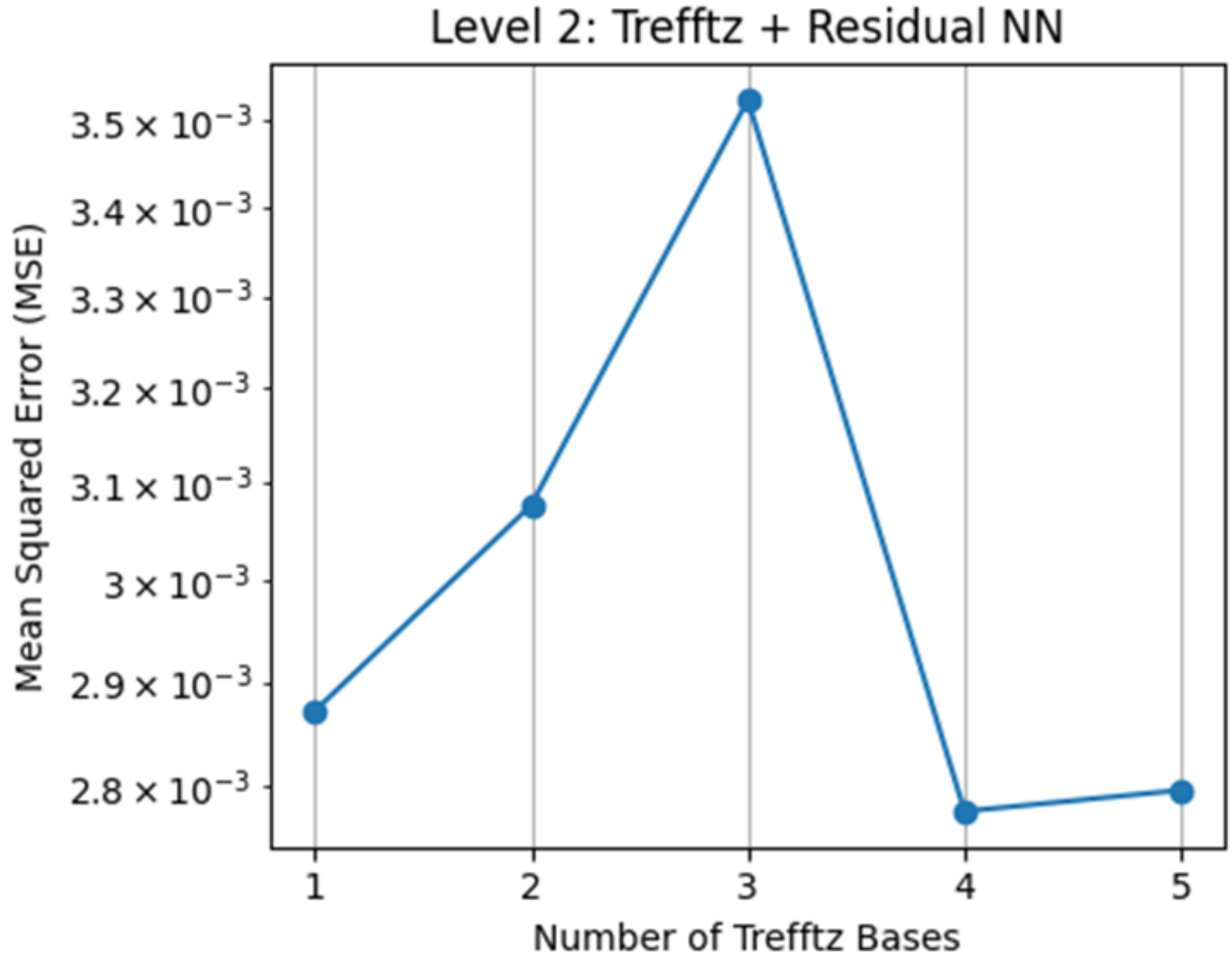

Figure 4: Relationship between the MSE and the number of Trefftz bases $\boldsymbol{N_b}$.

These results indicate that an appropriate balance between the number of Trefftz bases and the capacity of the residual network is essential. Importantly, even when MSE values are comparable across different configurations, the preservation of magnetic field-line structures remains superior in Trefftz-PINNs within the optimal basis range

### 4.3. CFD Streamline Comparison

The effectiveness of Trefftz-PINNs is further examined through a CFD surrogate modeling task using a decaying Taylor–Green vortex. Fig. 5 compares streamlines obtained from the exact solution, standard PINNs, and Trefftz-PINNs under matched MSE conditions.

Although standard PINNs achieve similar MSE levels, their reconstructed velocity fields exhibit phase shifts and distortions in streamline patterns. In some regions, streamlines deviate from their expected symmetric configuration, indicating a breakdown of coherent flow structures.

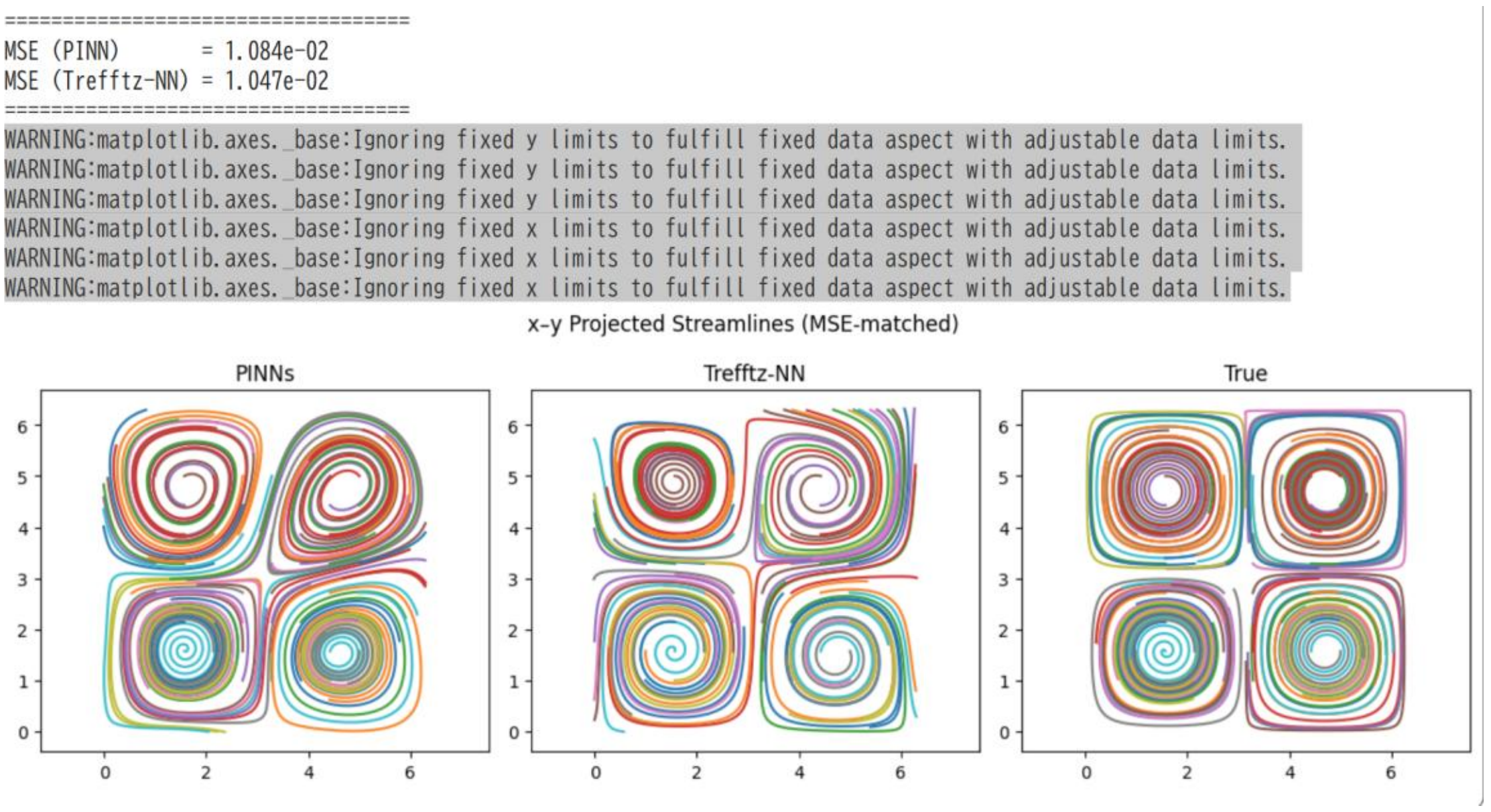

Figure 5: Comparison of streamlines obtained from the exact solution, standard PINNs, and Trefftz-PINNs under matched MSE conditions.

In contrast, Trefftz-PINNs reproduce streamlines that closely resemble those of the exact solution. The vortex symmetry and streamline continuity are well preserved, confirming that the physically admissible solution space imposed by Trefftz bases enhances structural robustness.

These results demonstrate that the advantages of Trefftz-PINNs are not limited to electromagnetic field problems but extend naturally to fluid flow analysis, where the preservation of flow topology is crucial.

## 5. Discussion

### 5.1. Why PINNs Fail Despite Small MSE

The present results clearly demonstrate that minimizing pointwise numerical error or PDE residuals alone does not guarantee the preservation of global physical structures. In both magnetic field and CFD examples, standard PINNs achieved MSE levels comparable to those of Trefftz-PINNs, yet exhibited severe structural breakdowns such as magnetic surface crossings and distorted streamlines.

This apparent contradiction arises from the nature of residual-based learning.

PINNs optimize a loss function defined over discrete sample points, typically combining PDE residuals and data fitting errors. However, these loss functions provide no explicit constraint on global topological properties of the solution. As a result, a neural network may converge to a function that satisfies the PDE approximately at sampled points while violating essential physical structures between them.

This phenomenon is closely related to the non-uniqueness of solutions under weak constraints. When boundary conditions, symmetries, or conservation properties are insufficiently enforced, the solution space contains many admissible functions with small residuals but fundamentally different global behaviors. PINNs, driven by gradient-based optimization, may converge to such spurious solutions without any explicit mechanism to reject them.

In this sense, the failure of PINNs observed in this study is not merely a numerical issue but a structural one. Even when residuals are small, the learned solution may lie outside the physically meaningful solution manifold.

## 5.2. Residual Hallucination as a Structural Failure Mode

The observed breakdowns can be interpreted as instances of *residual hallucination*, where a neural network produces a solution that appears valid according to residual-based metrics but is physically incorrect. Unlike classical numerical instability, residual hallucination originates from the mismatch between the optimization objective and the physical requirements of the problem.

In the heat conduction example presented in Section 1, the PINN converged to a solution with nearly zero Laplacian residual while deviating significantly from the exact temperature distribution. This illustrates that residual minimization alone is insufficient to recover the correct solution when gradients vanish or when the loss landscape contains flat regions. Similar effects are observed in the magnetic and CFD problems, where residual hallucination manifests as topological distortions rather than simple amplitude errors.

Importantly, residual hallucination is not an artifact of poor training or insufficient data. Even with identical training samples and matched MSE levels, the structural discrepancies persist. This indicates that residual hallucination is fundamentally linked to the expressiveness of the neural network and the absence of global physical constraints

## 5.3. Solution Space Restriction via Trefftz Bases

The Trefftz-PINN framework addresses this limitation by restricting the solution space *prior* to learning. By construction, Trefftz basis functions exactly satisfy the governing PDE, ensuring that the dominant component of the solution resides within the physically admissible solution manifold.

In this formulation, the neural network is not tasked with discovering the entire solution from scratch. Instead, it acts as a residual corrector that operates within a pre-constrained space. This significantly reduces the degrees of freedom available for generating unphysical solutions and suppresses residual hallucination at its source.

From a functional analysis perspective, this approach transforms the learning problem from an unconstrained approximation to a constrained one, where infeasible regions of the solution space are eliminated a priori. As a result, global structures such as magnetic surfaces and coherent vortices are preserved even when local numerical errors remain comparable.

This design philosophy contrasts with conventional PINNs, which rely on loss balancing, adaptive sampling, or optimization heuristics to guide the network toward physically meaningful solutions. Trefftz-PINNs instead enforce physical validity through basis selection, shifting the burden from training strategy to model architecture.

## 5.4. Relation to Existing Studies and Reduced-Order Modeling

Previous studies on PINN failures have primarily focused on improving training stability through adaptive loss weighting, curriculum learning, or modified network architectures. While these approaches address optimization challenges, they do not fundamentally restrict the solution space and therefore cannot fully prevent residual hallucination.

The present work complements these studies by highlighting the importance of solution space design. The Trefftz-PINN approach shares conceptual similarities with reduced-order modeling techniques such as Proper Orthogonal Decomposition (POD), where dominant physical modes are extracted to reduce dimensionality. However, unlike POD, Trefftz bases are derived from the governing equations themselves and therefore provide stronger physical guarantees.

Moreover, the Trefftz-PINN framework remains flexible and extensible. The residual neural network allows for modeling effects not captured by the selected bases, such as weak nonlinearity or modeling uncertainty, while preserving the core physical structure. This balance between physical rigor and learning flexibility is particularly advantageous for surrogate modeling in complex engineering systems.

### 5.5. Implications for Physics-Informed Machine Learning

The findings of this study suggest that future developments in physics-informed machine learning should place greater emphasis on *solution space construction* rather than solely on loss function design. Ensuring that learned solutions remain within physically admissible manifolds is essential for reliable surrogate modeling, especially when global structures play a critical role.

By demonstrating that Trefftz-PINNs preserve physical structures under matched MSE conditions, this work provides evidence that pre-constrained learning architectures offer a promising path toward robust and interpretable physics-informed models

## 6. Conclusion

This study demonstrated that achieving small numerical errors or PDE residuals does not necessarily guarantee the preservation of essential physical structures in physics-informed neural networks. Through matched-MSE comparisons, we showed that standard PINNs may produce physically implausible solutions, such as magnetic surface crossings and distorted streamlines, despite satisfying residual-based accuracy metrics.

To address this limitation, we introduced a Trefftz-PINN framework that constrains the solution space *prior* to learning by incorporating Trefftz basis functions that exactly satisfy the governing equations. Numerical experiments on vacuum magnetic field problems and CFD benchmark flows confirmed that Trefftz-PINNs consistently preserve global physical structures while maintaining comparable MSE levels to standard PINNs.

These results highlight residual hallucination as a fundamental failure mode of unconstrained residual minimization and emphasize the importance of solution space design in physics-informed machine learning. By shifting the focus from loss engineering to physically admissible representations, Trefftz-PINNs offer a robust and interpretable surrogate modeling paradigm.

The proposed framework is general and applicable beyond the test cases presented here, providing a promising foundation for reliable surrogate modeling in electromagnetic analysis, CFD, and complex multi-physics systems relevant to fusion reactor design.

## Acknowledgement

This work was performed with the support and under the auspices of the NIFS Collaboration

Research program (NIFS24KIEI001).

## References

[1] M. Raissi, P. Perdikaris, G. E. Karniadakis: Physics-informed neural networks: A deep learning framework for solving forward and inverse problems involving nonlinear partial differential equations, *Journal of Computational Physics*, 378:1 (2019), 686–707.

[2] M. Raissi, A. Yazdani, G. E. Karniadakis: Hidden physics models: Machine learning of nonlinear partial differential equations, *Journal of Computational Physics*, 357:15 (2018), 125–141.

[3] S. Wang, X. Yu, P. Perdikaris: When and why PINNs fail to train: A neural tangent kernel perspective, *Journal of Computational Physics*, 449:15 (2022), 110768.

[4] S. Wang, Y. Teng, P. Perdikaris: Understanding and mitigating gradient pathologies in physics-informed neural networks, *SIAM Journal on Scientific Computing*, 43:5 (2021), A3055–A3081.

[5] L. D. McClenny, U. Braga-Neto: Self-adaptive physics-informed neural networks using a soft attention mechanism, *Journal of Computational Physics*, 474:1 (2023), 111722.

[6] Z. Q. J. Xu, Y. Zhang, Y. Luo, Y. Xiao, Z. Ma: Frequency principle: Fourier analysis sheds light on deep neural networks, *Communications in Computational Physics*, 28:5 (2020), 1746–1767.

[7] S. Wang, S. Sankaran, P. Perdikaris: Respecting causality is all you need for training physics-informed neural networks, *arXiv*:2203.07404 (2022).

[8] A. S. Krishnapriyan, A. Gholami, S. Zhe, R. M. Kirby, M. W. Mahoney: Characterizing possible failure modes in physics-informed neural networks, in *Proc. NeurIPS*, Virtual Only, 2021.

[9] S. Mishra, R. Molinaro: Estimates on the generalization error of physics-informed neural networks, *ESAIM: Mathematical Modelling and Numerical Analysis*, 55:5 (2021), 1943–1965.

[10] G. E. Karniadakis, I. G. Kevrekidis, L. Lu, et al.: Physics-informed machine learning, *Nature Reviews Physics*, 3 (2021), 422–440.

[11] R. Akbarian Bafghi, et al.: Predictions and Corrections: Neural Predictors with Solver-Based Correction for ODEs and PDEs., in *Proceedings of NeurIPS (Workshop)*, 2025.

[12] E. Clement, et al: A Novel A.I Enhanced Reservoir Characterization with a

Combined Mixture of Experts -- NVIDIA Modulus based Physics Informed Neural Operator Forward Model*, arXiv*:2012.04426 (2024).

[13] P. Holmes, J. L. Lumley, G. Berkooz: *Turbulence, coherent structures, dynamical systems and symmetry*, Cambridge University Press, Cambridge, 1996.

[14] C. W. Rowley, S. T. M. Dawson: Model reduction for flow analysis and control, *Annual Review of Fluid Mechanics*, 49 (2017), 387–417.

[15] E. Trefftz: Ein Gegenstück zum Ritzschen Verfahren, in *Proc. the 2nd International Congress of Applied Mechanics*, Zürich, 1926, 131-137..

[16] Li, Z. C., & Lu, T. T. (2008). The Trefftz Method for Partial Differential Equations. CRC Press, ISBN 978-1584887539.

[17] R. Hiptmair, A. Moiola, I. Perugia: Trefftz methods for the Helmholtz equation, *Journal of Computational Physics*, 230:7 (2011), 2423–2440.

[18] L. Lu, R. Pestourie, W. Yao, Z. Wang, F. Verdugo, S. G. Johnson, Physics-informed neural networks with hard constraints for inverse design, *arXiv*:2102.04626 (2021).

[19] K. Koyamada, Y. Long, T. Kawamura, K. Konishi: Data-driven derivation of partial differential equations using neural network model, *International Journal of Modeling, Simulation, and Scientific Computing*, 12:2 (2021), 2140001.